\documentclass[a4paper,12pt]{article}
\usepackage[cp1251]{inputenc}
\usepackage[english, russian, ukrainian]{babel}
\usepackage{amsmath}
\usepackage{amssymb}

\newcommand{\1}{1\!\!\,{\rm I}}

\newcommand{\vf}{\varphi}

\newcommand{\mbR}{{\mathbb R}}
\newcommand{\mfX}{{\mathfrak X}}

\newcommand{\wt}{\widetilde}
\newcommand{\wh}{\widehat}
\renewcommand{\tilde}{\widetilde}
\renewcommand{\bar}{\overline}

\newcommand{\cE}{{\mathcal E}}

\newcommand{\pt}{\partial}
\newcommand{\Id}{\mathop{\rm Id}}

\begin{document}
 \large

\begin{center}
{\bf
Stochastic anticipating boundary value problems
}\\[1cm]
Andrey A. Dorogovtsev\\[1 cm]
Institute of Mathematics, Ukrainian Academy of Sciences,\\
ul. Tereshenkovskaia,3,Kiev, Ukraine\\
e-mail: adoro@imath.kiev.ua
\end{center}
\vskip1cm

{\bf Introduction.} This article is devoted to the stochastic
anticipating equations with the extended stochastic integral with
respect to the Gaussian processes of a special type. In the
particular cases the solutions of such an equations are the
well-known Wiener functionals after the second quantization. As an
application the stochastic Kolmogorov equation for the conditional
distributions of the diffusion process is obtained. Also we will
consider the conditional variant of the Feynman--Kac formula. The
two last sections of the article are devoted to the smoothing
problem in the case when noise is represented by the two jointly
Gaussian Wiener processes, which can have not a semimartingale
property with respect to the joint filtration. In order to describe
the objects of our consideration more explicitly consider the
following examples.

{\bf Example 0.1.}  Let $w_1, w_2$ be an independent one-dimensional
Wiener pro\-ces\-ses starting at zero. For fixed $x>0$ define
$$
\tau_x=inf\{t: x+w_2(t)=w_1(t)\}.
$$
If we rewrite $\tau_x$ in the form
$$
\tau_x=inf\{t: x=w_1(t)-w_2(t)\}
$$
then it is well-known [1] that $\tau_x$ is finite with probability
one. Denote by $u(x,t)=E\{\1_{\{\tau_x\leq t\}}/w_1\}.$ The
problem lies in the description of $u(x,t).$ One of the possible
approaches can be following. Denote
$$
\wt{w}(t)=\frac{1}{\sqrt{2}}(w_1(t)-w_2(t)).
$$
Then $\1_{\{\tau_x\leq t\}}$ is the functional from the
$\tilde{w}$ on $[0; t].$  So it has  an Ito-Wiener expansion as a
series of the multiple Wiener integrals:
$$
\1_{\{\tau_x\leq t\}}=
\sum^\infty_{k=0}
\mathop{\int\stackrel{k}{\ldots}\int}\limits_0^t
a_k(r_1,\ldots,r_k)d\tilde{w}(r_1)\ldots d\tilde{w}(r_k).
$$

Now, note that the conditional expectation of $\1_{\{\tau_x\leq t\}}$
with respect $w_1$ can be viewed as the action of the Ornstein-Uhlenbeck
semigroup operator $T_{\frac{1}{2}\ln2}$ [2] on $\1_{\{\tau_x\leq t\}}$.
Then
$$
u(x,t)
=\sum^\infty_{k=0}
\left(-\frac{1}{\sqrt{2}}\right)^k
\mathop{\int\stackrel{k}{\ldots}\int}\limits_0^t
a_k(r_1,\ldots,r_k)d{w}(r_1)\ldots d{w}(r_k).
$$
So the description of $u(x,t)$  can be obtained if we know the
Ito-Wiener expansion of $\1_{\{\tau_x\leq t\}}.$  One can receive
the kernels from this expansion in the following way. Consider for
the function $\vf\in L_2([0;t])$  the Wick exponent [2]
$$
e^{\int^t_0\vf(s)d\tilde{w}(s)-\frac{1}{2}\|\vf\|^2}=
\sum^\infty_{k=0}\frac{1}{k!}
\mathop{\int\stackrel{k}{\ldots}\int}\limits_0^t
\vf(r_1),\ldots,\vf(r_k)d\tilde{w}(r_1)\ldots d\tilde{w}(r_k).
$$
Here $\|\vf\|$  is the usual norm in $L_2([0;t]).$  It follows from
the properties of Ito-Wiener expansion, that
$$
E\1_{\{\tau_x\leq t\}}
e^{\int^t_0\vf(s)d\tilde{w}(s)-\frac{1}{2}\|\vf\|^2}=
\sum^\infty_{k=0}
\mathop{\int\stackrel{k}{\ldots}\int}\limits_0^t
\vf(r_1),\ldots,\vf(r_k)a_k(r_1,\ldots,r_k)dr_1\ldots dr_k.
$$
Consequently, for our purposes it is enough to know
$$
E\1_{\{\tau_x\leq t\}}
e^{\int^t_0\vf(s)d\tilde{w}(s)-\frac{1}{2}\|\vf\|^2}
$$
for the sufficiently large set of functions $\vf.$   It occurs that
the last  mathematical expectation can be obtained as a solution of
the certain boundary value problem for the parabolic PDE. This will
be done in the section 3 in more general situation.

{\bf Example 0.2}. Consider the ordinary stochastic differential
equation in $\mbR$
$$
dx(t)=a(x(t))dt+b(x(t))dw(t)
$$
with the smooth enough coefficients $a$ and $b.$  Denote by $x(r,s,t)$ the
solution which starts at the moment $s$ from the point $r.$ Let the function
$\vf\in C^2(\mbR)$  has bounded derivatives. Define for $r\in\mbR,
s\in[0;T]$
$$
\Phi(r,s)=\Gamma(\vf(r,s,T)),\eqno(0.1)
$$
where $\Gamma$ is the certain operator of the second quantization [2]. In
particular $\Gamma$  can be a mathematical expectation. Then, it can be
proved, that $\Phi$  satisfies the following partial stochastic differential
equation
$$
d\Phi(s,r)=-\left[
\frac{1}{2}b^2(r)\frac{\pt^2}{\pt r^2}\Phi(r,s)+a(r)\frac{\pt}{\pt r}\Phi(r,
s)\right] ds+
$$
$$
+b(r)\frac{\pt}{\pt r}\Phi(r,s)d\gamma(s).
$$
Here $\gamma(s)=\Gamma w(s)$  and the last differential is treated in the sennse
of anticipating stochastic integration. When $\Gamma$  is the mathematical
expectation, the last term vanishes.

This two examples shows the main goal of this article. Namely
there exist situations when the naturally arising Wiener
functionals satisfy the anticipating stochastic differential
equations and can be described with the using of stochastic
calculus. Here we propose the appropriate machinery and solve  the
correspondent equations.

We will consider the second quantization transformation of the
different Wiener functionals, whose mathematical expectations
usually arise in the probability representation of the solutions for
a boundary value problems. For such transformed functionals we will
get the anticipating stochastic equations with the extended
stochastic integral and solve it. Accordingly to this aim the
article is organized as follows. The first part contains the
properties of the second quantization operators in connection with
the extended stochastic integral or, more generally, with the
Gaussian strong random operators [3]. In the second part we consider
examples which show how the Fourier-Wiener transform works under
solution of the boundary value problems with the extended stochastic
integral on the interval. In the next part of the article we
introduce the anticipating boundary  value problems for some Wiener
functionals and discuss the solution of the such problems. The
section 4 and 5 are devoted to the following pair of equations
$$
\begin{aligned}
dx_1(t)=a_1(x_1(t))+dw_1(t)\\
dx_2(t)=a_2(x_1(t))+dw_2(t),
\end{aligned}
$$
where $w_1, w_2$ are jointly Gaussian Wiener processes, which can
have not a semimartingale property with respect to the joint
filtration. In this case we will look for the equation for
$E(f(x_1(t))/x_2).$

 {\bf 1. Second quantization and integrators.}   The material
of this section is partially based on the works [4,5]. Corresponding
facts are placed here for the completeness of the exposition but
their proofs are omitted. New claims are presented with the proofs.

We will start here with the abstract picture, when the ``white
noise'' generated by the Wiener process is substituted by the
generalized Gaussian random element in the Hilbert space. Let $H$
be a separable real Hilbert space with the   norm $\|\cdot\|$ and
inner product $(\cdot,\cdot).$  Suppose that $\xi$ is the
generalized Gaussian random element in $H$ with zero mean and
identical covariation. In other words $\xi$ is the family of
jointly Gaussian random variables denoted by $(\vf, \xi), \vf\in
H$ with the properties

1) $(\vf, \xi)$ has the normal distribution with zero mean and
variance $\|\vf\|^2$  for every $\vf\in H,$

2) $(\vf, \xi)$  is linear with respect to $\vf.$

During this section we suppose that all the random variables and
elements are measurable with respect to $\sigma(\xi)=\sigma((\vf,
\xi), \vf\in H.$ If the random variable $\alpha$ has the finite
second moment, than $\alpha$ has an Ito-Wiener expansion [6]
$$
\alpha=\sum^\infty_{k=0}A_k(\xi, \ldots,\xi). \eqno(1.1)
$$
Here, for every $k\geq1 \ A_k(\xi, \ldots, \xi)$   is the infinite-dimensional
generalization of the Hermite polinomial from $\xi,$   correspondent to the
$k$-linear symmetric Hilbert-Shmidt form $A_k$  on $H.$  Moreover, now the
following relation holds
$$
E\alpha^2=\sum^\infty_{k=0}k!\|A_k\|^2_k.
\eqno(1.2)
$$
Here $\|\cdot\|_k$  is the Hilbert-Shmidt form in $H^{\otimes k}.$  The same
expansion for $H$-valued random elements will be necessary. Let $x$  be a
random element in $H$  such, that
$$
E\|x\|^2<+\infty.
$$
Then, for every $\vf\in H$
$$
(x,\vf)=\sum^\infty_{k=0}A_k(\vf; \xi,\ldots,\xi).
\eqno(1.3)
$$
It can be easily checked using (1.2), that now $A_k$  is the
$k+1$-linear (not necessary symmetric) Hilbert-Shmidt form. So one
can write now
$$
x=\sum^\infty_{k=0}\tilde{A}_k(\xi,\ldots,\xi), \eqno(1.4)
$$
where
$$
\tilde{A}_k(\vf_1,\ldots,\vf_k):=
\sum^\infty_{j=1}A_k(e_j; \vf_1,\ldots, \vf_k)e_j
$$
for the arbitrary orthonormal basis $\{e_j; j\geq1\}$  in $H$   and the
series (1.4)  converges in $H$ in the square mean. The relation (1.2)
remains to be true
$$
E\|x\|^2=\sum^\infty_{k=0}k!\|\tilde{A}_k\|^2_k,
\eqno(1.5)
$$
where $\|\tilde{A}_k\|_k$  is the Hilbert--Shmidt norm in the space of
$H$-valued $k$-linear forms on $H.$

Now recall the definition of the operators of the second
quantization. Let $C$  be a continuous linear operator in $H.$
Suppose that the operator norm $\|C\|\leq1.$   Then for $\alpha$
and $x$  from (1.1)  and (1.4) define
$$
\begin{aligned}
&\Gamma(C)\alpha =\sum^\infty_{k=0}A_k(C\xi,\ldots,C\xi),\\
&\Gamma(C)x =\sum^\infty_{k=0}\tilde{A}_k(C\xi,\ldots,C\xi),\\
\end{aligned}
\eqno(1.6)
$$
where for $k\geq1 \ A_k(C\cdot,C\cdot,\ldots,C\cdot)$ and
$\tilde{A}_k(C\cdot,C\cdot,\ldots,C\cdot)$ are new Hilbert-Shmidt forms.

Using the estimation
$$
\|A_k(C\cdot,C\cdot,\ldots,C\cdot)\|_k\leq\|C\|^k\cdot\|A_k\|_k
$$
it is  easy to prove [2], that $\Gamma(C)$ is a continuous linear
operator in the space of square integrable random variables or
elements in $H.$

{\bf Definition 1.1.} [2] Operator $\Gamma(C)$  is the operator of \
t\,h\,e \ s\,e\,c\,o\,n\,d \
q\,u\-\,a\,n\,\-t\,i\,\-z\,a\,\-t\,i\,o\,n \ correspondent to the
operator $C.$

 Before to consider some examples, we will present the useful
representation of the  second quantization operators. Let $\xi'$  be
the generalized Gaussian random element in $H,$  independent and
equidistributed with $\xi.$

Consider the following generalized Gaussian random element in $H$
$$
\eta=\sqrt{1-CC^*}\xi'+C\xi.
\eqno(1.7)
$$
This element can be properly defined by the formula
$$
\forall \vf\in H: \ \ \ \ \
(\vf, \eta):=(\sqrt{1-CC^*}\vf,\xi')+(C^*\vf,\xi).
\eqno(1.8)
$$
Note that $\eta$  has zero mean and identity covariation. In order to check
this,  is sufficient to note the relation
$$
\|\vf\|^2=\|\sqrt{1-CC^*}\vf\|^2+\|C^*\vf\|^2.
$$
For every random variable $\alpha$ with an expansion (1.1)  define
$$
\alpha(\eta):=\sum^\infty_{k=0}A_k(\eta,\ldots,\eta).
$$
The following representation will be useful.

{\bf Lemma 1.1.} {\sl For arbitrary $\alpha\in L_2(\Omega, \sigma(\xi), P)$
and operator $C$ in $H$ with $\|C\|\leq1$
$$
\Gamma(C)\alpha =E(\alpha(\eta)/\xi).
\eqno(1.9)
$$
}

{\it Proof}. Note, that the both parts of (1.9) are continuous
with respect to  $\alpha$ in the square mean. So, it is enough to check
(1.9)  for the following random variables
$$
e^{(\vf,\xi)-\frac{1}{2}\|\vf\|^2}, \ \vf\in H. \eqno(1.10)
$$
Really,  the random variable of this kind has the Ito-Wiener
expansion of the form
$$
e^{(\vf,\xi)-\frac{1}{2}\|\vf\|^2}=\sum^\infty_{k=0}\frac{1}{k!}
\vf^{\otimes k}(\xi,\ldots,\xi).
$$
Here $\vf^{\otimes k}$ is $k$-th tensor power of $\vf$  which acts on
$H$ by the rule
$$
\vf^{\otimes k}(\psi_1,\ldots,\psi_k)=\prod^k_{j=1}(\vf,\psi_j).
$$
So, as it was mentioned in introduction, for $\alpha,$ which has an
expansion (1.1)
$$
E\alpha e^{(\vf,\xi)-\frac{1}{2}\|\vf\|^2}=\sum^\infty_{k=0}A_k
(\vf,\ldots,\vf).
$$
Hence, the set of all linear combinations of the variables (1.10)
is dense in $L_2.$ Now [2] the following equality holds
$$
\Gamma(C)e^{(\vf,\xi)-\frac{1}{2}\|\vf\|^2}=
e^{(C^*\vf,\xi)-\frac{1}{2}\|C^*\vf\|^2}. \eqno(1.11)
$$
From other side
$$
e^{(\vf,\eta)-\frac{1}{2}\|\vf\|^2}=
e^{(C^*\vf,\xi)-\frac{1}{2}\|C^*\vf\|^2}\cdot
e^{(\sqrt{1-CC^*}\vf,\xi')-\frac{1}{2}\|\sqrt{1-CC^*}\vf\|^2}.
$$
In order to finish the proof it is enough now note, that
$$
E
e^{(\sqrt{1-CC^*}\vf,\xi')-\frac{1}{2}\|\sqrt{1-CC^*}\vf\|^2}=1,
$$
and that $\xi'$ and $\xi$  are independent. It follows from here, that
$$
E\left(e^{(\vf,\eta)-\frac{1}{2}\|\vf\|^2}/\xi\right)=
e^{(C^*\vf,\xi)-\frac{1}{2}\|C^*\vf\|^2}.
$$
Lemma is proved.

This lemma has the following useful application for us.

{\bf Corollary 1.1.} {\sl
Let $\Gamma(C)$ be an operator of the second quantization. Let $x$  be $a$
random element in the complete separable metric space $\mfX$
measurable with respect to $\xi.$   Then there exists the random
probability measure $\mu$ on $\mfX$  such, that for every bounded
measurable function $f: \mfX\to\mbR$   the following equality holds
$$
\int_\mfX fd\mu =\Gamma(C)f(x).
$$
}

{\it Proof.} Let us define $\mu$  as a conditional distribution of $x(\eta)$
with respect to $\xi.$  Then, for every bounded measurable $f: \mfX\to\mbR$
$$
\int_\mfX fd\mu=E(f(x(\eta))/\xi)=\Gamma(C)f(x).
$$
The unique difficulty on this way lies in the proper definition of
$x(\eta)$ (remind, that $\xi$ and $\eta$ are not usual random
elements). In order to break this difficulty we will use the
following analog of the Levy theorem. Let $\{e_j; j\geq1\}$   be
an orthonormal basis in $H.$  Define the sequences of random
elements in $H$  by the rule
$$
\begin{aligned}
&
\xi_n=\sum^n_{j=1}(e_j,\xi)e_j,\\
&
\eta_n=\sum^n_{j=1}(e_j,\eta)e_j, \ n\geq1.
\end{aligned}
$$
Note that the sequences $\{\xi_n;n\geq1\}$  and $\{\eta_n;n\geq1\}$ are
equidistributed. Now for every $n\geq1$ consider the random measure
$\nu_n$ in $\mfX,$  which is built in the following way
$$
\nu_n(\Delta)=E\{\1_\Delta(x)/\xi_n\}.
$$
Here $\Delta$ is an arbitrary Borel subset of $\mfX.$  This random
measures have two important properties. First of all for every
$n\geq1$ $\nu_n$ can be viewed as $\widetilde{\nu}_n(\xi_n),$
where $\widetilde{\nu}_n$ is a Borel function from $H$ to the
space of all probability measures on $\mfX$ equipped with the
distance of weak convergence. Secondly, with probability one
$\nu_n$ weakly converge to $\delta_x$ under $n$ tends to infinity.
The last assertion follows from the usual Levy theorem [7]. More
precisely, for arbitrary continuous bounded function $f:
\mfX\to\mbR$
$$
f(x)=E(f(x)/\xi)=\lim_{n\to\infty}E(f(x)/\xi_n)=
$$
$$
=\lim_{n\to\infty}\int_\mfX f(u)\nu_n(du) \ \mbox{a.s.}
$$
Taking $f$  from the countable set which define the weak
convergence [8] we get the required statement. Now note, that the
sequence of random measures $\{\widetilde{\nu}_n(\eta\; n\geq1\}$
is equidistributed with $\{\widetilde{\nu}_n(\xi_n); n\geq1\}.$
Hence with probability one there exists the weak limit of
$\widetilde{\nu}_n(\eta_n)$  which is a delta-measure concentrated
in the certain random point $y.$  This random point $y$ is by
definition $x(\eta).$  The correctness of this definition can be
easily checked. Lemma is proved.

Consider the examples of the random measures, which arise in the
application of the Corollary 1.1 and will be important for us.

{\bf Example 1.1.}  Suppose, that $H=L_2([0;T], \mbR^d).$  Define the
generalized Gaussian random element $\xi$ in $H$  with the help of the
$d$-dimensional Wiener process $W$ on $[0; T].$  Namely, for
$
\vf=(\vf_1,\ldots,\vf_d)\in L_2([0;1], \mbR^d)$ define
$$
(\vf,\xi):=\sum^d_{j=1}\int^T_0\vf_j(s)dW_j(s).
\eqno(1.12)
$$
Now consider the functions $a: \mbR^d\to\mbR^d$ and $b:
\mbR^d\to\mbR^{d\times d}$ which satisfy the Lipshitz condition
and the domain $G$ in $\mbR^d$  with the $C^1$-boundary $\Gamma$.
Let for every $s\in[0;T]$  and $u\in\mbR^d$ $x(u,s,T)$  denote the
solution at time $T$ of the following Cauchy problem
$$
\begin{cases}
dx(t)=a(x(t))dt+b(x(t))dW(t),\\
x(s)=u.
\end{cases}
\eqno(1.13)
$$
Denote by $\nu_{u,s}$  the random measure obtained from $x(u,s,T)$
via corollary 1.1 with the help of the certain operator of the
second quantization $\Gamma(C).$  In the next section we will obtain
the stochastic variant of the Kolmogorov equation for $\nu_{u,s}.$
Note, that in the case $C=0$ measures $\nu_{u,s}$ became to be
deterministic and satisfy the usual Kolmogorov equation.

Now let us define for every $u\in G$  the random moment
$$
\tau_{u,s}=\inf\{T, t\leq T: x(u,s,t)\in\Gamma\}.
$$
Let $\mu_{u,s}$ be the random measure obtained from
$x(u,s,\tau_{u,s})$ via corollary 1.1. It occurs, that measures
$\mu_{u,s}$ satisfy certain anticipating boundary  value problem.

In order to describe the anticipating SPDE for the random measures
from above mentioned example we need in the relation between the
operators of the second quantization and extended stochastic
integral. We will study this connection in the more general
situation when the extended stochastic integral is substituted by
the general Gaussian strong random operator (GSRO in the sequel).
Let us recall the following definition.

{\bf Definition 1.2}. [3] T\,h\,e \ G\,a\,u\,s\,s\,i\,a\,n \
s\,t\,r\,o\,n\,g \ r\,a\,n\,d\,o\,m \ l\,i\,n\,e\,a\,r \
o\,p\,e\,\-r\,a\,\-t\,o\,r \ (\,G\,S\,R\,O\,) \ $A$ in $H$ is the
mapping, which maps every element $x$ of $H$  into the jointly
Gaussian with $\xi$ random element in $H$ and is continuous in the
square mean.

As an example of GSRO the integral with respect to Wiener process can be
considered.

{\bf Example 1.2.} Consider $H$  and $\xi$  from example 1.1. Let for
simplicity $d=1.$  Define GSRO $A$ in the following way
$$
\forall \vf\in H: \ \ \ \ \ \
(A\vf)(t)=\int^t_0\vf(s)dw(s), \ \ t\in[0; T].
$$
It can be easily seen that $A\vf$ now is a Gaussian random element in
$H,$  and $A$ is continuous in square mean.

For to include in this picture the  integration with respect to
another Gaussian processes (for example with respect to the
fractional Brownian motion) consider more general GSRO. Suppose,
that $K$ be a bounded linear operator, which acts from
$L_2([0;T])$  to $L_2([0;T]^2).$   Define
$$
\forall \vf\in H: \ \ \ \
(A\vf)(t)=\int^T_0(K\vf)(t,s)dw(s).
$$
It can be checked, that $A$ is GSRO in $H.$  Making an obvious changes one can
define the GSRO acting from the different Hilbert space $H_1$  into $H.$
For example consider for $\alpha\in\left(
\frac{1}{2};1\right)$ the covariation function of the fractional Brownian
motion [9]  with Hurst parameter $\alpha$
$$
R(s,t)=\frac{1}{2}(t^{2\alpha}+s^{2\alpha}-|t-s|^{2\alpha}).
$$
Define the space $H_1$  as a completion of the set of step
functions on $[0;T]$ with respect the inner product under which
$$
(\1_{[0;s]},\1_{[0;t]})=R(s,t).
$$
Consider the kernel $K^\alpha$ from the integral representation of the
fractional
Brownian motion $B^\alpha$ [10]
$$
B^\alpha(t)=\int^t_0 K^\alpha(t,s)dw(s)
$$
and
$$
\frac{\pt K^\alpha}{\pt t}(t,s)=c_\alpha\left(\alpha-\frac{1}{2}\right)
(t-s)^\alpha-\frac{3}{2}\left(\frac{s}{t}\right)^{\frac{1}{2}-\alpha}.
$$
Define for $\vf\in H_1$
$$
(K\vf)(t,s)=\int^t_s\vf(r)\frac{\pt K^\alpha}{\pt r}(r,s)dr\1_{[0;t]}(s).
$$
Now let
$$
(A\vf)(t)=\int^T_0(K\vf)(t,s)dw(s)=\int^t_0(K\vf)(t,s)dw(s).
$$
Then
$$
(A\vf)(t)=\int^t_0\vf(s)dB^\alpha(s).
$$

We will consider the action of GSRO on the random elements in $H.$
Corresponding definition was proposed in [3, 11]. Consider
arbitrary GSRO $A$ in $H.$ Then for every $\vf\in H$ the
Ito-Wiener expansion of $A\vf$ contains only two terms
$$
A\vf=\alpha_0\vf+\alpha_1(\vf)(\xi).\eqno(1.14)
$$
Here $\alpha_0$  is a continuous linear operator in $H$  and $\alpha_1$ is
a continuous linear operator from $H$ to the space of Hilbert-Shmidt
operators in $H.$ Now let $x$  be a random element in $H$  with the finite
second moment. Then $\alpha_1(x)$ has a finite second moment in the space
of Hilbert-Shmidt operators. So for every $\vf\in H$
$$
\alpha_1(x)(\vf)=\sum^\infty_{k=0}B_k(\vf; \xi,\ldots,\xi).
$$
It can be easily verified, that $B_k$ is $k+1$-linear $H$-valued
Hilbert-Shmidt form on $H.$ Define $\Lambda B_k$ as a symmetrization of
$B_k$  with respect to all $k+1$ variables.

{\bf Definition 1.3.} [3, 11]  The random element $x$ \
b\,e\,l\,o\,n\,g\,s \ t\,o \ t\,h\,e \ d\,o\,m\,a\,i\,n \ o\,f \
d\,e\,f\,i\,n\,i\,t\,i\,o\,n \ o\,f \ G\,S\,R\,O \ $A$ if the series
$$
\sum^\infty_{k=0}\Lambda B_k(\xi,\ldots,\xi)
$$
converges in $H$  in the square mean and in this case
$$
Ax=\alpha_0x+\sum^\infty_{k=0}\Lambda B_k(\xi,\ldots,\xi).
\eqno(1.15)
$$

In the partial cases this definition gives us the definition of
the extended stochastic integral [6, 12, 13]. We will define this
integral for the special class of Gaussian processes. Suppose,
that $H=L_2([0;T])$  and $\xi$ is generated by the Wiener process
$w$ as above. Consider the jointly Gaussian with $w$  process
$\{\gamma(t); t\in[0;T]\}$ with zero mean.

{\bf Definition 1.4.} [14]  Process $\gamma$ is \ a\,n \
i\,n\,t\,e\,g\,r\,a\,t\,o\,r \ if there exists the constant $C$ such
that for every step function $\vf$ on $[0;T]$
$$
\vf=\sum^{n-1}_{k=0}a_k\1_{[t_k; t_{k+1})}
$$
the following unequality holds
$$
E(\sum^{n-1}_{k=0}a_k(\gamma(t_{k+1})-\gamma(t_k)))^2\leq C
\sum^{n-1}_{k=0}a_k^2(t_{k+1}-t_k).
\eqno(1.16)
$$

The good examples of the integrators can be obtained via the following
simple  statement.

{\bf Lemma 1.2.} {\sl Let $\Gamma(C)$ be an operator of the second
quantization. Then
$$
\gamma(t)=\Gamma(C)w(t), t\in[0;T]
$$
is an integrator.
}

The proof of this lemma easily follows from the properties of $\Gamma(C).$

Note that the integrator can have the unbounded quadratic
variation and consequently can have not the semimartingale
properties (see [14]). It is easy to see from (1.16), that for
every integrator $\gamma$ and $\vf\in L_2([0;T])$ the stochastic
integral
$$
\int^t_0\vf d\gamma
$$
exists as a limit of the integrals from the step functions and
$$
E\left(\int^t_0\vf d\gamma\right)^2\leq C\int^t_0\vf^2(s)ds.
$$
So one can define GSRO $A_\gamma$  associated with the integrator $\gamma$
by the rule
$$
\forall \vf\in L_2([0;T]): \ \ \ \ \
(A_\gamma\vf)(t)=\int^t_0\vf d\gamma.
$$
In this situation the definition 1.3 became to be a definition of the
extended stochastic integral with respect to $\gamma.$  Note that in the
case $\gamma=w$ it will be a usual extended integral.

Now let us consider the relation between the action of GSRO and the
operators of the second quantization.

{\bf Theorem 1.1.} [4]  {\sl  Let $A$ be a GSRO in $H$ and
$\Gamma(C)$ be an operator of the second quantization. Suppose
that the random element $x$  lies in the domain of definition of
$A$ in the sence of definition 1.3. Then $\Gamma(C)x$ belongs to
the domain of definition of GSRO $\Gamma(C)A$  and the following
equality holds
$$
\Gamma(C)(Ax)=\Gamma(C)A (\Gamma(C)x).
\eqno(1.17)
$$
Here $\Gamma(C)A$ is the GSRO which acts by the rule
$$
\forall\vf\in H: \ \Gamma(C)A\vf=\Gamma(C)(A\vf).
$$
}

The proof of this theorem is placed in [4] and so is omitted.
Instead the proof consider the following important example of
application of the theorem 1.1 to the stochastic integration.

{\bf Example 1.3.}  Consider in the situation of the example 1.2
GSRO of integration with respect to Wiener process $w.$   Suppose
that random function $x$ in $L_2([0;T])$  with the finite second
moment is adapted to the flow of $\sigma$-fields generated by $w.$
It is well-known [12, 13], that in this case the extended
stochastic integral
$$
\int^t_0x(s)dw(s), \ t\in[0;T]
$$
exists and coincides with the  Ito integral. Now the theorem 1.1
says us that
$$
\Gamma(C)\left(\int^t_0x(s)dw(s)\right)=\int^t_0\Gamma(C)x(s)d\gamma(s),
$$
where $\gamma$ is an integrator of the type
$\gamma(t)=\Gamma(C)w(t)$  and the integral in the right part is
an extended stochastic integral.

{\bf 2. Anticipating equations and boundary value problems for the
second quantization of the certain Wiener functionals.}   In this
section we use the theorem 1.1 for to get the stochastic equations
for the functions from example 1.1. In order to do this we will
need in the notion of the weak solution of the equation with the
extended stochastic integral. Let us consider for GSRO $A$, random
element $x$  from the domain of $A$  and $\vf\in H$
$$
E(Ax)e^{(\vf,\xi)-\frac{1}{2}\|\vf\|^2}.
\eqno(2.1)
$$
This expectation in Hilbert space is a Bochner integral.

{\bf Lemma 2.1.} {\sl
Integral (2.1) can be expressed in the following way
$$
E(Ax)e^{(\vf,\xi)-\frac{1}{2}\|\vf\|^2}=
\alpha_0x_\vf+\alpha_1(x_\vf)(\vf).
\eqno(2.2)
$$

Here $\alpha_0$ and $\alpha_1$ are the terms from expansion of $A$ (1.14) and
$$
x_\vf=Exe^{(\vf,\xi)-\frac{1}{2}\|\vf\|^2}.
$$
}

{\it Proof.}  Consider the Ito-Wiener expansion of $x$
$$
x=\sum^\infty_{k=0}B_k(\xi,\ldots,\xi).
$$
It is known [2], that
$$
e^{(\vf,\xi)-\frac{1}{2}\|\vf\|^2}=
\sum^\infty_{k=0}\frac{1}{k!}\vf^{\otimes k}(\xi,\ldots,\xi).
$$
So, by the properties of Ito-Wiener expansion
$$
x_\vf=\sum^\infty_{k=0}B_k(\vf,\ldots,\vf).
$$
By the usual continuity arguments it is enough to verify the statement of
the lemma for the case $\alpha_0=0$   and $x=B_k(\xi,\ldots,\xi)$
for a certain $k\geq0.$  In this case
$$
Ax=\Lambda D_{k+1}(\xi,\ldots,\xi),
$$
where the $k+1$-linear Hilbert-Shmidt form $D_{k+1}$ has a following
representation
$$
D_{k+1}(\vf_1,\ldots,\vf_{k+1})=\alpha_1(B_k(\vf_1,\ldots,\vf_k)) (\vf_{k+1}).
$$
Since
$$
D_{k+1}(\vf,\ldots,\vf)=\alpha_1(x_\vf)(\vf)
$$
the lemma is proved.

This lemma allows to define the action of the unbounded Gaussian
random operator on the random elements in the weak sense. Let
$H_1$ be a separable real Hilbert space densely and continuously
embedded into $H.$ Consider GSRO $A$ acting from $H_1$ to $H$ (in
another word for every $u\in H_1$ $Au$ is a Gaussian element in
$H$ and this correspondence is continuous in the square mean with
respect to the convergence in $H_1$). Then $A$  can be treated as
an unbounded Gaussian random operator in $H.$  As it was mentioned
above, $A$ can be described with the help of two deterministic
linear operators $\alpha_0: H_1\to H,$ $\alpha_1:
H_1\to\sigma_2(H)$   (here $\sigma_2(H)$ is the space of the
Hilbert-Shmidt operators in $H$ with the correspondent norm).
Operators $\alpha_0$ and $\alpha_1$  can be considered as an
operators on $H.$  Then the correspondence
$$
H\supset H_1\ni u\mapsto\alpha_0(u)+\alpha_1(u)(\xi)
$$
defines an unbounded Gaussian random operator in $H$  with the domain of
definition $H_1.$  The next definition is closely related to the
lemma 2.1.

{\bf Definition 2.1.}   The random element $x$ in $H$ with the
finite second moment \ b\,e\,l\,o\,n\,g\,s \ i\,n  \ t\,h\,e  \
w\,e\,a\,k \ s\,e\,n\,s\,e \ t\,o \ t\,h\,e \ d\,o\,m\,a\,i\,n \
o\,f \ d\,e\,f\,i\,n\,i\,t\,i\,o\,n \ o\,f \ t\,h\,e \
u\,n\,b\,o\,u\,n\,d\,e\,d \ G\,a\,u\,s\,s\,i\,a\,n \
r\,a\,n\,d\,o\,m \ o\,p\,e\,r\,a\,t\,o\,r \ $A$ if there exist the
dense subset $L\subset H$ and the random element $y$ in $H$ with the
finite second moment such, that
$$
\forall\vf\in L: \ \ \
x_\vf=Exe^{(\vf,\xi)-\frac{1}{2}\|\xi\|^2}\in H_1, \
\alpha_0(x_\vf)+\alpha_1(x_\vf)(\vf)=y_\vf. \eqno(2.3)
$$
Here $y_\vf$  is defined in the same way as $x_\vf.$

Consider an example of the situation described in the previous
definition. Note, that in a similar way the definition of the action of
unbounded Gaussian random operator on the Hilbert space $H'\ne H$  can be
defined.

{\bf Example 2.1.}  Case of deterministic operator. Let
$H'=L_2(\mbR\times[0;1], e^{-\frac{x^2}{2}}dx\times dt).$
Consider the nonrandom operator $A$ in $H'$  which is defined by the
formula
$$
Af(x,t)=\int^t_0\frac{\pt}{\pt x}f(x,s)ds.
$$
As a Hilbert space $H_1$ let us use
$$
H_1=W^1_2(\mbR, e^{-\frac{x^2}{2}}dx)\times L_2([0;1]),
$$
i.e. $H_1$  consists of the functions from $x\in\mbR, t\in[0;1]$   which
have one Sobolev derivative with respect to $x$ which is in $H'.$
Consider the following random element in $H'$
$$
X(x,t)=\1_{\{w(t)\leq x\}}.
$$
It is easy to verify that $X$  has a finite second moment
$$
E\int_{\mbR}\int^1_0X(x,t)^2\cdot e^{-\frac{x^2}{2}}dxdt\leq\sqrt{2\pi}.
$$
Note, that for fixed $t$  $X(\cdot, t)\notin
W^1_2(\mbR,e^{-\frac{x^2}{2}}dx).$ Indeed due to the Sobolev
embedding theorem all functions from $W^1_2(\mbR,
e^{-\frac{x^2}{2}}dx)$ must be continuous. Let us prove that $X$
belongs to the domain of definition $A$ in the sense of definition
2.1. Take $\vf\in L_2([0;1])$ and consider
$$
X_\vf(x,t)=E\1_{\{w(t)\leq x\}}\exp\left\{
\int^1_0\vf(s)dw(s)-\frac{1}{2}\int^1_0\vf^2(s)ds
\right\}.
$$
It follows from Girsanov theorem that
$$
X_\vf(x,t)=E\1_{\{w(t)+\int^t_0\vf(s)ds\leq x\}}.
$$
Then $X_\vf\in H_1$ and
$$
\int^t_0\frac{\pt}{\pt x}X_\vf(x,s)ds=\int^t_0\frac{1}{\sqrt{2\pi
s}} \exp-\frac{1}{2s} \left\{x-\int^s_0\vf(r)dr\right\}^2ds.
\eqno(2.4)
$$
Now find $Y$ with the finite second moment in $\wt{H}$  such, that
$Y_\vf$  equals to (2.4). Consider the sequence
$$
Y^n(x,t)=2n\int^t_0\1_{[x-\frac{1}{n};x+\frac{1}{n}]}(w(s))ds,
n\geq1.
$$
Note [1], that there exists the random field $Y(x,t), x\in\mbR,
t\in[0;1]$ such, that
$$
E(Y^n(x,t)-Y(x,t))^2\to0, n\to\infty.
$$
Also note, that
$$
\sup_{n\geq1} EY^n(x,t)^2\leq c|x|.
$$
This can be easily checked using Tanaka approach. Hence $Y$ is the random
element in $\wt{H}$ with the finite second moment. Moreover
$$
Y_\vf(x,t)=EY(x,t) \exp \left\{
\int^1_0\vf(s)dw(s)-\frac{1}{2}\int^t_0\vf^2(s)ds \right\}=
$$
$$
=\lim_{n\to\infty }EY^n(x,t)
\exp
\left\{
\int^1_0\vf(s)dw(s)-\frac{1}{2}\int^t_0\vf^2(s)ds
\right\}=
$$
$$
=\lim_{n\to\infty}E2n
\int^t_0\1_{[x-\frac{1}{n};x+\frac{1}{n}]}(w(s))ds
\exp
\left\{
\int^1_0\vf(s)dw(s)-\frac{1}{2}\int^t_0\vf^2(s)ds
\right\}=
$$
$$
=\lim_{n\to\infty}E2n
\int^t_0\1_{[x-\frac{1}{n};x+\frac{1}{n}]}(w(s)+\int^s_0\vf(r)dr)ds=
$$
$$
=\int^t_0\frac{1}{\sqrt{2\pi s}}\exp-\frac{1}{2s}
\left\{
x-\int^s_0\vf(r)dr
\right\}^2ds.
$$
So,
$$
Y_\vf(x,t)=\int^t_0\frac{\pt}{\pt x}X_\vf(x,s)ds.
$$
Consequently $X$  belongs to the domain of definition of operator
$A$  in the sense of definition 2.1  and $AX$  is the local time of
the Wiener process. Note, that the Ito-Wiener expansion of the
Wiener local time and  more general related objects can be found in
[15, 16].

{\bf Example 2.2}.  Let the spaces $\wt{H}$  be as in the previous example.
Suppose, that the space $H_1$  is built similarly to the previous
example with the substitution of $W^1_2(\mbR, e^{-\frac{x^2}{2}}dx)$  by
$W^2_2(\mbR, e^{-\frac{x^2}{2}}dx).$  Define the Gaussian random operator
$A$ on $H_1$ in the following way
$$
Af(x,t)=\frac{1}{2}\int^t_0\frac{\pt^2}{\pt x^2}f(x,s)ds+\int^t_0
\frac{\pt}{\pt  x}f(x,s)dw(s).
$$
Let us consider arbitrary bounded and measurable function $h$ on $\mbR$
and
define the random element $X$ in $\wt{H}$  by the formula
$$
X(x,t)=h(x+w(t)).
$$
Prove, that $X$  belongs to the domain of definition of $A$ in the weak
sense. Really, the operators $\alpha_0$  and $\alpha_1$  from the Ito-Wiener
representation of $A$ now has the form
$$
\alpha_0(u)(x,t)=\frac{1}{2}\int^t_0\frac{\pt^2}{\pt x^2}u(x,s)ds,
$$
$$
\alpha_1(u)(\vf)(x,t)=\int^t_0\frac{\pt}{\pt x}u(x,s)\vf(s)ds.
$$
Consider
$$
X_\vf(x,t)=
Eh(x+w(t))\exp
\left\{
\int^1_0\vf(s)dw(s)-\frac{1}{2}\int^t_0\vf^2(s)ds\right\}=
$$
$$
=
Eh(x+w(t)+\int^t_0\vf(s)ds)=
\frac{1}{\sqrt{2\pi t}}\int_\mbR h(u)
\exp-\frac{1}{2t}
\left\{
u-x-\int^t_0\vf(s)ds
\right\}^2du.
$$
Consequently,
$$
\int^t_0\frac{1}{2}\frac{\pt^2}{\pt x^2}X_\vf(x,s)+
\frac{\pt}{\pt x}X_\vf(x,s)\vf(s)ds =X_\vf(x,t)-h(x).
$$
The last equality is valid every sure with respect to the Lebesque
measure. So it is right as an equality between the elements from
$\wt{H}.$ Finally we see, that $X$  belongs to the domain of
definition $A$  in the weak sense.

Now let us turn to the equations with the GSRO. Here we will
consider two types of anticipating equations. In the first one the
equation has the usual Ito character but the initial value depends
on future increments of Wiener process.  In the second one the
equation itself contains the integral with respect to the
integrators.

Firstly let us consider the equation with the operator from the example 2.2
written in the differential form
$$
\begin{cases}
dU(x,t)=\frac{1}{2}\frac{\pt^2}{\pt x^2}U(x,t)dt+\frac{\pt}{\pt x}U(x,t)
dw(t),\\
U(x,0)=f(x), x\in\mbR.
\end{cases}
\eqno(2.5)
$$
Suppose that initial condition $f$ is the random function, which is
measurable with respect to the Wiener process
$\{w(t), t\in[0;1]\}.$  The equation is treated in the weak sense. It is
well-known, that the linear SDE with anticipating initial condition
$$
\begin{cases}
dz(t)=az(t)dt+bz(t)dw(t),\\
z(0)=\alpha,
\end{cases}
\eqno(2.6)
$$
where the stochastic integral is treated in the Skorokhod sense has the
solution
$$
z(t)=Q_t\alpha\cdot\cE^t_0.
$$
Here $\cE^t_0$ is the usual stochastic exponent related to (2.6)
and
$$
Q_t\alpha=\alpha(w(\cdot)+b\cdot\wedge t).
$$
I.e. $Q_t$  is the shift on the Wiener space. In the dimension
greater then one  analog of this statement is not valid in general
[17]. The reason is that the matrix-valued coefficients $a$ and
$b$  can be noncommuting and the direction of the shift can not be
obtained. But in this case another approach can be used [4]. This
approach is related to the von-Neimann series for the solution of
the equation with GSRO. We will consider this approach for
equation (2.5). Suppose, that the initial condition  $f$  in (2.5)
has the form
$$
f(x)=\alpha\cdot\vf(x),
$$
where $\vf$  is the deterministic function from $L_2(\mbR)$  and
$\alpha$ is the random variable measurable with respect to $w.$
Consider the Fourier transform $\wh{U}$ of the solution. Then it
satisfies following Cauchy problem
$$
\begin{cases}
d\wh{U}(\lambda,t)=-\frac{1}{2}\lambda^2\wh{U}(\lambda,t)dt+i\lambda
\wh{U}(\lambda,t)dw(t),\\
\wh{U}(\lambda,0)=\alpha\wh{\vf}(\lambda).
\end{cases}
\eqno(2.7)
$$
This problem can be treated under fixed $\lambda.$  Let us write the
correspondent equation
$$
\begin{cases}
dz(t)=az(t)dt+ibz(t)dw(t),\\
z(0)=\alpha.
\end{cases}
\eqno(2.8)
$$
This equation differs from (2.6)  due to the presence of $i.$  Consider the
real and imaginary parts of $z$ separately.

Let $z(t)=z_1(t)+iz_2(t),$  where $z_1,z_2$ are the real stochastic
processes. Then
$$
\begin{cases}
dz_1(t)=az_1(t)dt-bz_2(t)dw(t),\\
dz_2(t)=az_2(t)dt+bz_1(t)dw(t),\\
z_1(0)=\alpha, \ z_2(0)=0.
\end{cases}
\eqno(2.9)
$$
Last equation can be written in the matrix form. Denote
$\vec{z}(t)=(z_1(t), z_2(t)).$ Then
$$
d\vec{z}(t)=a\Id\vec{z}(t)dt+B\vec{z}(t)dw(t),
\eqno(2.10)
$$
where $\Id$ is an identical $2\times2$ matrix and
$$
B=\begin{pmatrix}
0&-b\\
b&0
\end{pmatrix}.
$$
Since $B$  has not now the system of the real eigen-values, then we
can not apply the formula with the shift of probability argument
from [17]. Let the matrix-valued process $Y(s,t), 0\leq s\leq t\leq
T$  be the solution of the Cauchy problem
$$
\begin{cases}
dY(s,t)=aY(s,t)dt+BY(s,t)dw(t)\\
Y(s,s)=\Id.
\end{cases}
$$
Suppose, that $\alpha$ has a finite Ito-Wiener expansion, i.e.
$$
\alpha=\sum^N_{k=0}
\mathop{\int\stackrel{k}{\ldots}\int}\limits_0^T
a_k(t_1,\ldots,t_k)dw(t_1)\ldots dw(t_k).
$$
The next statement was proved in [4].

{\bf Theorem 2.1} [4]. {\sl The Cauchy problem (2.9) has the
unique solution of the form
$$
\vec{z}(t)=
\sum^\infty_{k=0}
\mathop{\int\stackrel{k}{\ldots}\int}\limits_{0\leq t_1
\leq\ldots\leq t_k\leq t}
D^k\alpha(\tau_1,\ldots,\tau_k)Y(t_k,t)B\cdot
$$
$$
\cdot
Y(t_{k-1},t_k)B\ldots BY(0,t_1)\vec{u}dt_1\ldots dt_k.
\eqno(2.11)
$$
Here for $k=0$ the  corresponding summand is $\alpha
Y(0,t)\vec{u}$   and $\vec{u}=(1,0).$ Since in our situation
$Y(s,t)$ has an evolutionary property, i.e. for $t_1\leq t_2\leq
t_3$
$$
Y(t_1,t_3)=Y(t_2,t_3)Y(t_1,t_2)
$$
and $B$ commute with $Y,$  then (2.11)  can be rewritten in the form
$$
\vec{z}(t)=
\sum^\infty_{k=0}Y(0,t)B^k\vec{u}
\mathop{\int\stackrel{k}{\ldots}\int}\limits_{0\leq t_1
\leq\ldots\leq t_k\leq t}
D^k\alpha(t_1,\ldots,t_k)
dt_1\ldots dt_k.
\eqno(2.12)
$$
}
Now let us recall, that vector $\vec{z}(t)$  consists of the real and
imaginary part of the solution to the initial Cauchy problem (2.8). In
order to obtain $z(t)$  consider the vector
$\vec{y}_k(t)=Y(0,t)B^k\vec{u}, k\geq1.$  In can be easily checked, that
$$
y_k(t)_1+iy_k(t)_2=i^kb^ke^{at+ibw(t)+\frac{1}{2}b^2t}.
$$
So
$$
z(t)=e^{at+\frac{1}{2}b^2t+ibw(t)}
\sum^\infty_{k=0}
i^kb^k
\mathop{\int\stackrel{k}{\ldots}\int}\limits_{0\leq t_1
\leq\ldots\leq t_k\leq t}
D^k\alpha(t_1,\ldots,t_k)
dt_1\ldots dt_k.
\eqno(2.13)
$$
Since the action of GSRO on the random element is closed in the
sense  of the linear operator theory, then (2.13)  remains valid
for those $\alpha,$  who has infinitely many stochastic
derivatives and for which the series converges in the square mean.
Returning to the equation (2.7) we obtain the explicit expression
of $\wh{U}(\lambda,t)$
$$
\wh{U}(\lambda,t)=\wh{\vf}(\lambda) e^{i\lambda w(t)}
\sum^\infty_{k=0}
i^k\lambda^k
\mathop{\int\stackrel{k}{\ldots}\int}\limits_{0\leq t_1
\leq\ldots\leq t_k\leq t}
D^k\alpha(t_1,\ldots,t_k)
dt_1\ldots dt_k.
\eqno(2.14)
$$
Taking the inverse Fourier transform
$$
U(x,t)=
\sum^\infty_{k=0}
\vf^{(k)}(x+w(t))
\mathop{\int\stackrel{k}{\ldots}\int}\limits_{0\leq t_1
\leq\ldots\leq t_k\leq t}
D^k\alpha(t_1,\ldots,t_k)
dt_1\ldots dt_k.
\eqno(2.14')
$$
This is    the solution of the initial anticipating SPDE, which is obtained
using GSRO approach.

Now let us consider the anticipating SPDE with the integrators.
Let $\{w(t); t\geq0\}$  be a Brownian motion in $\mbR^d,$  $G$  is
a domain in $\mbR^d$   with the smooth boundary. Consider the
diffusion process $\{x(t); t\geq0\}$ in $G$  with the diffusion
matrix $(a_{ij}(x))^d_{ij=1},$  drift term $(b_i(x))_{i=1,d}.$ We
will consider two types of the conditions on the boundary. First
is the reflection. In this case we can suppose that $x$ satisfies
the Skorokhod SDE in $G$
$$
\begin{cases}
dx(t)=b(x(t))dt+\sigma(x(t))dw(t)+\1_{\pt G}(x(t))\gamma(x(t))d\eta(t),\\
x(0)=x_0, \eta(0)=0.
\end{cases}
\eqno(2.15)
$$
Here $\sigma(x)\sigma^*(x)=(a_{ij}(x)),$  $\gamma$ is the direction of the
reflection and $\eta$ is the local time of $x$ on the boundary
$\pt G.$  The second type of the process is the diffusion inside $G$  with
the coefficients $a$  and $b$ stopped on the boundary. Let us denote this
process by $y.$  We will suppose that $y$ satisfies the same SDE as $x$
from (2.15)  with out the last term. Consider for the function
$f\in C^2(\bar{G})$  and the certain bounded linear operator $C$  in
$L_2([0;T], \mbR^d)$   following random function
$$
U(x_0,t)=\Gamma(G)f(x(t)).
$$
The function $U$ is defined on $G\times[0;+\infty).$  Of course we will
suppose, that $\|C\|\leq1$ in order to define $U$ correctly.

{\bf Theorem 2.2.} {\sl Suppose, that $\sigma, b$ and $\gamma$
have three bounded continuous derivatives and boundary of $G$ is
the $C^3$-bounded manifold. Let also $(\gamma, \nabla f)|_{\pt
G}=0.$  Then $U$ satisfies  the following anticipating PSDE
$$
\begin{gathered}
 dU(x,t)= \left[\frac{1}{2}\sum^d_{ij=1}a_{ij}(x)\frac{\pt^2}{\pt
x_i\pt x_j} U(x,t)+\sum^d_{j=1}b_j(x)\frac{\pt}{\pt
x_j}U(x,t)\right]dt-\\
 -\sum^d_{ij=1}\sigma_{ij}(x)\frac{\pt}{\pt
x_i}U(x,t)d\gamma_j(t), \end{gathered}
\eqno(2.16)
$$
and the boundary condition
$$
(\gamma, \nabla U)|_{\pt G}=0.
$$
Here for $j=1,\ldots,d$
$$
\gamma_j(t)=\Gamma(G)w_j(t).
$$
The equation (2.16)  is understood in the weak sense, i.e. the
random function $U$  as a random element in the space
$L_2(G\times[0;T])$ belongs to the domain of definition of the
random integral operator from the right part of (2.16) via the
definition 2.1. }

{\it Proof.}  Consider $\vf\in C^1([0;T], \mbR^d).$  Denote
$$
\cE^t_0(\vf)=
\exp\left\{
\int^T_0\sum^d_{j=1}\vf_j(s)dw_j(s)-
\frac{1}{2}\int^T_0\sum^d_{j=1}\vf^2_j(s)ds\right\}.
$$
Now, due to (1.11)
$$
E\Gamma(C)f(x(t))E^t_0(\vf)=
$$
$$
=Ef(x(t)) \exp\left\{ \int^T_0\sum^d_{j=1}(C\vf)_j(s)dw_j(s)-
\frac{1}{2}\int^T_0\sum^d_{j=1}(C\vf)^2_j(s)ds\right\}.
$$

Let us denote for the convenience $C\vf=\psi.$  Note, that
$\cE^T_0(\psi)$  is the probability density from Girsanov theorem.
So,
$$
Ef(x(t))\cE^T_0(\psi)=
Ef(x(t))\cE^t_0(\psi)=Ef(z(t)),
$$
where the process $z$ satisfies the following Cauchy problem
$$
\begin{cases}
dz(t)=[b(z(t))+\sigma(z(t))\psi(t)]dt+\\
+\sigma(z(t))dw(t)+\1_{\pt G}(z(t))\gamma(z(t))d\zeta(t),\\
z(0)=x_0, \zeta(0)=0,
\end{cases}
$$
$\zeta$ is the local time of $z$ on the boundary. Denote
$$
V(x_0,t)=Ef(z(t)).
$$

The function $V$  satisfies the boundary condition
$(\gamma, \nabla V)|_{\pt G}=0$  and the partial differential equation
$$
\frac{\pt}{\pt t}V(x,t)=\frac{1}{2}\sum^d_{ij=1}a_{ij}(x)
\frac{\pt^2}{\pt x_i\pt x_j}V(x,t)+
$$
$$
+\sum^d_{j=1}b_j(x)\frac{\pt}{\pt x_j}V(x,t)+
\sum^d_{ij=1}\sigma_{ij}(x)\psi_j(t)\frac{\pt}{\pt x_i}V(x,t).
\eqno(2.17)
$$

Note, that the right part of (2.17) is exactly the expression from
definition 2.1. Theorem is proved.

{\bf Remark.}  Due to the uniqueness of the solution to (2.17) the solution
of the initial anticipating SPDE also is unique.

Consider now the functionals from the diffusion process $y.$  Denote by $\tau$
the hitting time of $y$ on the boundary $\pt G.$ Let $g\in C(\pt G).$
Consider the following functional
$$
\int^\tau_0f(y(s))ds+g(y(\tau)).
$$
Our aim is to get the equation for the random function
$$
\Gamma(C)\left(
\int^\tau_0f(y(s))ds+g(y(\tau))
\right).
\eqno(2.18)
$$
Unfortunately we can not obtain the equation for (2.18)  itself.
Instead we will consider the following function
$$
Q(x,t)=\Gamma(C)
\left(
\int^{T\wedge\tau}_t f(y(s))ds+\wt{g}(y(T\wedge\tau))
\right).
\eqno(2.19)
$$
Here $\wt{g}$  is the unique deterministic function from
$C^2(G)\bigcap C(\bar{G})$  which satisfies the Dirichlet problem
$$
\begin{cases}
\frac{1}{2}\sum^d_{ij=1}a_{ij}(x)\frac{\pt^2}{\pt x_i\pt x_j}\wt{g}(x)+
\sum^d_{j=1}b_j(x)\frac{\pt\wt{g}(x)}{\pt x_j}=0, x\in G,\\
\wt{g}|_{\pt G}=g.
\end{cases}
\eqno(2.20)
$$
The following statement holds

{\bf Theorem 2.3.} {\sl The random function $Q$  satisfies in the
weak sense the following anticipating SPDE with the boundary
conditions
$$
\begin{cases}
dQ(x,t)=
\left[
-\frac{1}{2}\sum^d_{ij=1}a_{ij}(x)\frac{\pt^2}{\pt x_i\pt x_j}Q(x,t)-
\sum^d_{j=1}b_j(x)\frac{\pt}{\pt x_j}Q(x,t)+f(x)\right]dt-\\
-\frac{1}{2}\sum^d_{ij=1}a_{ij}(x)\frac{\pt}{\pt x_j}Q(x,t)d\gamma_i(t),\\
Q|_{\pt G}=g, \ Q(x,T)=\wt{g}(x).
\end{cases}
\eqno(2.21)
$$
Here $\gamma_i, i=1,\ldots,d$   are the same as in the previous theorem and
the  stochastic integrals related to $d\gamma_i$  have the same meaning.
}

The proof of this theorem is almost analogous to the proof of the
theorem 2.2 and is omitted.

{\bf 3. Second quantization of the hitting moment}. Consider the application of
the results from previous section to the example from the introduction.
Let $\{w(t); t\geq0\}$  be the one-dimensional Wiener process starting from
the point $x>0.$  Denote by $\tau$ the hitting time for $w$ on the level 0.
Let $\Gamma(C)$  be the certain operator of the second quantization. We are
interesting in the value $\Gamma(C)\1_{\{\tau\leq t\}}.$  Firstly let us
find $\Gamma(C)f(w(t\wedge \tau)),$  where
$f\in C^2([0;+\infty])$  bounded and satisfies condition $f''(0)=0.$   Let us
denote
$$
\xi_{x,s}(t)=x+w(t)-w(s), \ 0\leq s\leq t,
$$
$$
\tau_{x,s}=\inf\{t\geq s: \xi_{x,s}(t)=0\},
$$
$$
V(x,s)=\Gamma(C)f(\xi_{x,s}(t\wedge\tau_{x,s})).
$$
Similar to theorem 2.1  $V$ satisfies the following anticipating SPDE
$$
\begin{cases}
dV(x,s)=-\frac{1}{2}\frac{\pt^2}{\pt x^2}V(x,s)ds-\frac{\pt}{\pt x}V(x,s)
d\gamma(s),\\
\frac{\pt^2}{\pt x^2}V(0,s)=0, \ V(x,t)=f(x).
\end{cases}
\eqno(3.1)
$$
Consider some special case operator $C.$

{\bf Example 3.1}.   Let $C$  be a projector on the certain unit vector $e$
in $L_2([0;T]), C=e\otimes e.$   Then $\Gamma(C)$  is the conditional
mathematical expectation with respect to the random variable
$$
\eta=\int^T_0 e(s)dw(s).
$$
So
$$
\gamma(t)=M(w(t)/\eta)=\eta\int^t_0e(s)ds.
$$
In this case the extended stochastic integral with respect to $\gamma(t)$
from the random process $\zeta,$  which has the stochastic derivative, has
the form
$$
\int^T_0\zeta(s)d\gamma(s)=\int^T_0(\zeta(s)\eta-(D\zeta(s),e))e(s)ds.
$$
Note, that in (3.1) $V$ is measurable with respect to $\gamma,$ i.e. (in
our case) with respect to $\eta.$  Let us try to find $V$ as a smooth
function not only $x$ and $s$  but also $\eta.$   Then (3.1) can be
rewritten as follows
$$
\begin{cases}
dV(x,s,\eta)=
\left[
-\frac{1}{2}\frac{\pt^2}{\pt x^2}V(x,s,\eta)-\frac{\pt}{\pt x}(\eta V(x,s,
\eta) -\frac{\pt}{\pt\eta}V(x,s,\eta))\right]ds, \\
\frac{\pt^2}{\pt x^2}V(0,s,\eta)=0, \ V(x,t,\eta)=f(x).
\end{cases}
\eqno(3.2)
$$
Substituting
$$
\wt{V}(x,s,\eta)=e^{-\frac{\eta^2}{2}}V(x,s,\eta)
$$
we get
$$
\begin{cases}
d\wt{V}(x,s,\eta)=
\left[
-\frac{1}{2}\frac{\pt^2}{\pt x^2}\wt{V}(x,s,\eta)
+\frac{\pt^2}{\pt x\pt\eta}\wt{V}(x,s,
\eta)\right]ds, \\
\frac{\pt^2}{\pt x^2}\wt{V}(0,s,\eta)=0, \ \wt{V}(x,t,\eta)=
e^{-\frac{\eta^2}{2}}f(x).
\end{cases}
\eqno(3.3)
$$

The last equation can bee solved by using the Ito-Wiener expansion of
$\wt{V}.$  Let us look for the solution of the kind
$$
\wt{V}(x,s,\eta)=\sum^\infty_{k=0}e^{-\frac{\eta^2}{2}}H_k(\eta)V_k(x,s).
\eqno(3.4)
$$
Here $H_k$ is the Hermite polynomial of order $k$  with the first
coefficient equal to one:
$$
H_k(\alpha)=(-1)^ke^{\frac{\alpha^2}{2}}\left(\frac{d}{d\alpha}\right)^k
e^{-\frac{\alpha^2}{2}}.
$$
Using the well-known relation
$$
\frac{d}{d\eta}e^{-\frac{\eta^2}{2}}H_k(\eta)=-e^{-\frac{\eta^2}{2}}
H_{k+1}(\eta)
$$
one can get the following system of the boundary value problems for the
deterministic coefficients in (3.4)
$$
\begin{cases}
\frac{\pt}{\pt s}V_0(x,s)=-\frac{1}{2}\frac{\pt}{\pt x^2}V_0(x,s),\\
\frac{\pt^2}{\pt x^2}V_0(0,s)=0, \
V_0(x,t)=f(x),
\end{cases}
\eqno(3.5)
$$
$$
\begin{cases}
\frac{\pt}{\pt s}V_{k+1}(x,s)=-\frac{1}{2}\frac{\pt}{\pt x^2}
V_{k+1}(x,s)-\frac{\pt}{\pt x}V_k(x,s),\\
\frac{\pt^2}{\pt x^2}V_{k+1}(0,s)=0, \
V_{k+1}(x,t)=0, k\geq0.
\end{cases}
\eqno(3.6)
$$
The obtained system can be solved recurrently. To do this denote for
$x,y\geq0$
$$
q_t(x,y)=\frac{1}{\sqrt{2\pi t}}
\left(
e^{-\frac{(x-y)^2}{2t}}-e^{-\frac{(x+y)^2}{2t}}\right).
$$
$q_t$  is the transition density for killed at zero Brownian motion. Now
$$
V_0(x,s)=\int^{+\infty}_0q_{t-s}(x,y)f(y)dy,
$$
$$
V_{k+1}(x,s)=\int^t_s\int^{+\infty}_0
q_{t-\tau}(x,y)\frac{\pt}{\pt y}V_k(y,\tau)dyd\tau, k\geq0.
$$
From this relations we get
$$
V_k(x,s)=
\mathop{\int\stackrel{k}{\ldots}\int}\limits_
{s\leq\tau_1\leq\ldots\leq\tau_k\leq t}
\mathop{\int\stackrel{k}{\ldots}\int}\limits_0^{+\infty}
q_{t-\tau_k}(x,y_1)
\frac{\pt}{\pt y_1}q_{\tau_k-\tau_{k-1}}(y_1,y_2)\ldots
$$
$$
\ldots\frac{\pt}{\pt y_{k-1}}q_{\tau_{k-\zeta}}
(y_{k-1},y_k)f(y_k)dy_1\ldots dy_kd\tau_1\ldots d\tau_k.
$$
Now the solution of (3.2) and (3.3)  can be written as a series
(3.4).

{\bf4. Smoothing problem.}  The last sections of the article are
devoted to the following problem. Let $(w_1,w_2)$ be the pair of
jointly Gaussian one-dimensional Wiener processes. Let the processes
$x_1, x_2$ are obtained via the relations
$$
\begin{array}{l}
dx_1(t)=a_1(x_1(t))dt+dw_1(t),\\
dx_2(t)=a_2(x_1(t))dt+dw_2(t),\\
x_1(0)=x_2(0)=0.
\end{array}
\eqno(4.1)
$$

Note, that the second equality is just a definition of $x_2$ but not
an equation. The problem is to find the conditional distribution of
$x_1(t)$ for $t\in[0;1]$ under given $\{x_2(s); s\in[0;1]\}.$
 We will try to get the equation for
 $$
 E(f(x_1(t))/x_2)
 $$
 for the appropriate functions $f.$

 Firstly let us study the joint distribution of $(w_1,w_2)$. Note,
 that there exists the bounded linear operator $V: L_2([0;1])\to
 L_2([0;1])$ such, that \newline
 $\forall\vf_1,\vf_2\in L_2([0;1]):$
 $$
 E\int^1_0\vf_1dw_1\int^1_0\vf_2dw_2=\int^1_0\vf_1 V\vf_2 ds.
 $$

 This fact follows from the reason, that the left part of the above
 formula is the continuous bilinear form with respect to $\vf_1$ and
 $\vf_2.$  Moreover, the operator norm $\|V\|\leq 1.$

  In
this section we consider the density of the distribution $(x_1,
x_2)$  with respect to the distribution of $(w_1,w_2)$ and study its
properties under the conditional expectation. The problem is that
the distribution of $(w_1,w_2)$ is not a Wiener measure in $C([0;1],
\mbR^2).$  So in order to get the density we need to adapt the
general Gaussian measure setup [18] to the our case. For  the future
let us denote $C([0;1])$ as $C$  and identify the space $C([0;1],
\mbR^2)$  with the direct sum $C\oplus C,$ which is furnished by the
sum of the norms. Denote also by $H$  the space
$$
L_2([0;1],\mbR^2)=L_2([0;1])\oplus L_2([0;1])
$$
with the scalar product defined by the formula
$$
(\vf, \psi)=\int^1_0\vf_1\psi_1ds+\int^1_0\vf_2\psi_2ds.
$$
With the pair $(w_1,w_2)$  we can associate the generalized Gaussian
random element $\xi$ in $H$  by the rule
$$
(\vf, \xi)=\int^1_0\vf_1dw_1+\int^1_0\vf_2dw_2.
$$
Note that $\xi$  has not an identity covariation operator. Really
$$
E(\vf,\xi)(\psi, \xi)=\int^1_0\vf_1\psi_1ds+\int^1_0\vf_2\psi_2ds+
$$
$$
+\int^1_0\vf_1 V\psi_2ds +\int^1_0\psi_1 V\vf_2ds.
$$
Here $V$  is described above bounded linear operator in
$L_2([0;1]).$ Denote by $S$ the operator in $H$  which acts by the
rule
$$
S\vf=(\vf_1+V\vf_2, V^*\vf_1+\vf_2).
$$
Then
$$
E(\vf, \xi)(\psi,\xi)=(S\vf, \psi).
$$
Our aim is to describe the transformations of the pair $(w_1,w_2)$
in the terms of $\xi.$  Let us start with the deterministic
admissible shifts. Denote by $i$ the canonical embedding of $H$ into
$C^2,$ i.e.
$$
i(\vf)(t)=(\int^t_0\vf_1ds,\int^t_0\vf_2ds).
$$

{\bf Lemma 4.1.} {\it Let the operator norm $\|V\|<1,$  then for
every $h\in H$  $i(h)$ is admissible shift for $\mu_{w_1,w_2}$  and
the corresponding density has the form
$$
p_h(\xi)=\exp\{(S^{-1}h,\xi)-\frac{1}{2}(S^{-1}h,h)\}. \eqno(4.2)
$$
}

{\bf Remark.}   Due to the condition $\|V\|<1$ the operator $S^{-1}$
is bounded on $H$  and can be written in the form
$$
S^{-1}\vf=\vf+(Q_{11}\vf_1+Q_{12}\vf_2,
Q_{21}\vf_1+Q_{22}\vf_2)=\vf+Q\vf,
$$
where $\|Q\|<1.$

{\it Proof.}  Note, that by the definition the operator $S$ is
nonnegative. Define
$$
{\xi}'=S^{-\frac{1}{2}}\xi.
$$
Then the shift of the distribution of $(w_1,w_2)$ on the vector
$i(h)$  is related to the shift of ${\xi}'$  on the vector
$S^{-\frac{1}{2}}h.$   Now the statement of the lemma follows from
the well-known formula for the density in the terms of ${\xi}'$
$$
p({\xi}')=\exp\{({\xi}',
S^{-\frac{1}{2}}h)-\frac{1}{2}(S^{-\frac{1}{2}}h,
S^{-\frac{1}{2}}h)\}
$$
if we rewrite it in the terms of $\xi.$

The lemma is proved.

{\bf Remark.}  Formula (4.2) can be rewritten in the terms of
$w_1,w_2.$ Really, by the definition
$$
(S^{-1}h,\xi)-\frac{1}{2}(S^{-1}h,h)=
$$
$$
= \int^1_0(S^{-1}h)_1dw_1+ \int^1_0(S^{-1}h)_2dw_2- \frac{1}{2}
\int^1_0(S^{-1}h)_1h_1ds- \frac{1}{2} \int^1_0(S^{-1}h)_2h_2ds=
\eqno(4.3)
$$
$$
= \int^1_0(h_1+Q_{11}h_1+Q_{12}h_2)dw_1+
\int^1_0(h_2+Q_{21}h_1+Q_{22}h_2)dw_2-
$$
$$
-\frac{1}{2} \int^1_0(h_1+Q_{11}h_1+Q_{12}h_2)h_1ds- \frac{1}{2}
\int^1_0(h_2+Q_{21}h_1+Q_{22}h_2)h_2ds.
$$

Using the same method one can find the density of $\mu_{x_1,x_2}$
with respect  $\mu_{w_1,w_2}.$  Firstly define the stochastic
derivatives of the functionals from $w_1,w_2$  with respect to $\xi$
and the extended stochastic integral in the terms of $\xi.$  Let
$\vf$  be a differentiable bounded function on $C\oplus C.$   Define
the stochastic derivative of the random variable $\vf(w_1,w_2)$  by
the formula
$$
D\vf(w_1,w_2):=i^*\nabla\vf(w_1,w_2).
$$
By this definition for every $t\in[0;1]$
$$
Dw_1(t)=(\1_{[0;t]},0),
$$
$$
Dw_2(t)=(0, \1_{[0;t]}).
$$
Note, that $\vf(w_1,w_2)$   can be regarded as a functional from the
generalized random element ${\xi}'$  which was introduced in the
proof of the lemma 4.1. Since ${\xi}'$  has an identity covariation
operator the stochastic derivatives and extended stochastic integral
for the functionals from ${\xi}'$  are connected by the usual
relation. Now we will define the extended stochastic integral with
respect to $\xi.$  It can be done in the following way. Consider the
Gaussian random functional on $H$  of the kind
$$
J(\vf)=(\vf, \xi).
$$
Then, in the terms of ${\xi}'$  $J$  can be rewritten as
$$
J(\vf)=(S^{\frac{1}{2}}\vf, {\xi}').
$$
So, the action of $J$  on the random element $x$  in $H$  via the
definition 1.3 has the form
$$
J(x)=I(S^{\frac{1}{2}}x). \eqno(4.4)
$$
Here $I$ is the extended stochastic integral with respect to
${\xi}'.$ Note also, that for the stochastic derivatives with
respect to ${\xi}'$  and ${\xi}$ we have the obvious relation
$$
D_{{\xi}'}\alpha=S^{-\frac{1}{2}}D_\xi\alpha.
$$
Hence on the domain of definition
$$
E(D_\xi\alpha, x)=E(S^{-\frac{1}{2}}D_\xi\alpha, S^{\frac{1}{2}}x)=
$$
$$
= E(D_{{\xi}'}\alpha, S^{\frac{1}{2}}x) =E\alpha\cdot
I(S^{\frac{1}{2}}x)= E\alpha J(x). \eqno(4.5)
$$
Thus the relation between the stochastic derivative and extended
stochastic integral with respect to $\xi$ is the same as for
${\xi}'.$ Now let us turn to the nonadapted shifts of the
distribution of $(w_1,w_2).$ Consider the pair of random processes
$x_1,x_2$  which are defined by the equations (4.1).

 The next lemma is standard.

{\bf Lemma 4.2.} {\it Let the functions $a_1,a_2$  be continuously
differentiable and have bounded derivatives. Then

1) for every $t\in[0;1]$  the random variables $x_1(t), x_2(t)$
have the stochastic derivatives $Dx_1(t), Dx_2(t),$

2) the random element $(a_1(x_1(\cdot)), a_2(x_1(\cdot)))$ in $H$
has the stochastic derivative, and
$$
D(a_1(x_1(s)), a_2(x_1(s)))(t)=
$$
$$
= (a'_1(x_1(s))Dx_1(s)(t), a'_2(x_1(s))Dx_1(s)(t)),
$$

3) the stochastic derivative of $x_1$  with respect to $w_1$  (i.e.
the first coordinate of $Dx_1$)  satisfies the equation
$$
\begin{array}{l}
D_1x_1(s)(t)=1+\int^s_ta'_1(x_1(r))D_1x_1(r)(t)dr, 0\leq t\leq s\leq1,\\
D_1x_1(s)(t)=0, \ \ t>s,
\end{array}
$$
and
$$
Dx_1(s)(t)=(D_1x_1(s)(t),0).
$$
}

It follows from the lemma 4.2, that $\|D(a_1(x_1(\cdot)),
a_2(x_1(\cdot)))\|_H$ can be made small if we take $a_1'$ and $a_2'$
small enough. Since the operator $S^{-\frac{1}{2}}$  is bounded in
$H,$  then due to the theorem 3.2.2 from [ 18 ]  the distribution of
$(x_1,x_2)$  is absolutely continuous with respect to the
distribution $(w_1,w_2)$   for sufficiently small $a_1',$ $a_2'.$
The corresponding density will be denoted by $p.$ Accordingly to
[18]  $p$  has the form
$$
p=\zeta\cdot\exp\{I(S^{-\frac{1}{2}}h)-\frac{1}{2}(S^{-1}h,h)\},
\eqno(4.6)
$$
where
$$
h(t)=(a_1(w_1(t)), a_2(w_1(t))),
$$
and $\zeta$  is the corresponding Carleman--Fredgolm determinant.
Due to (4.4) (4.6) can be rewritten as
$$
p=\zeta\exp\{J(S^{-1}h)-\frac{1}{2}(S^{-1}h, h)\}. \eqno(4.7)
$$
This expression allows us to conclude, that up to the term $\zeta$
$p$  has the stochastic derivative. We will suppose, that this is so
in the next section, where the formulas for the conditional
expectation and extended stochastic integral will be obtained in
non-Gaussian case. \vskip15pt

{\bf5. Conditional expectation.}  For the processes $(x_1,x_2)$ from
(4.1) let us search for the conditional distribution of $x_1(t)$
under fixed $\{x_2(s); s\in[0;1]\}.$   Firstly note, that under our
conditions the distribution of $(x_1,x_2)$  is absolutely continuous
with respect to the distribution of $(w_1,w_2)$  and, consequently,
the distribution  of $x_2$  is absolutely continuous with respect to
the distribution of $w_2.$

Denote for a moment by $\mu$  the distribution of the pair
$(w_1,w_2)$  on $C([0;1])\oplus C([0;1])$   and by $\mu_1$  and
$\mu_2$  the distributions of $w_1$  and $w_2$  (shurely these are a
Wiener measures but on the different  copies of $C([0;1]).$ It
follows from the general theory of integration   that the measure
$\mu$   can be desintegrated with respect to $\mu_2,$  i.e.
$$
\mu(\Delta)=\int_{C([0;1])}\nu(u,\Delta_u)\mu_2(du),
$$
for arbitrary Borel $\Delta$  in $C([0;1])\oplus C([0;1]).$  Here
$\nu$  is a measurable family of the probability measures and
$\Delta_u=\{v\in C([0;1]): (v,u)\in \Delta\}.$

Define for the measurable bounded function $\vf: C([0;1])\to\mbR$
the function $\psi$  by the rule
$$
C([0;1])\ni u\mapsto\psi(u)=\int_{C([0;1])}\vf(v)p(v,u)\nu(u,dv).
\eqno(5.1)
$$
$$
\cdot (\int_{C([0;1])}p(v,u)\nu(u, dv))^{-1}.
$$

The following variant of the Bayes formula holds.

{\bf Lemma 5.1.}
$$
E(\vf(x_1)/x_2)=\psi(x_2).
$$

{\it Proof.}   Note firstly, that $\psi(x_2)$  is correctly defined
because the function $\psi$  is defined up to the set of Wiener
measure zero, and the distribution of $x_2$ is equivalent to this
measure. Now for arbitrary bounded and  measurable function $\gamma:
C([0;1])\to\mbR$
$$
E\vf(x_1)\gamma(x_2)=E\vf(w_1)\gamma(w_2)p(w_1,w_2)=
$$
$$
=E\gamma(w_2)\cdot E(\vf(w_1)p(w_1,w_2)/w_2)=
$$
$$
 =E\gamma(w_2)\cdot \frac
{E(\vf(w_1)p(w_1,w_2)/w_2)}{E(p(w_1,w_2)/w_2)}\cdot
E(p(w_1,w_2)/w_2)=
$$
$$
= E\gamma(w_2)\psi(w_2)\cdot p(w_1,w_2)=E\gamma(x_2)\psi(x_2).
$$
This finishes the proof.

For arbitrary $t\in[0;1]$  denote by $\pi_t$   the random measure on
$\mbR$   whose pairing with the bounded measurable function $f$  is
defined by the formula
$$
\int_{\mbR} f(r)\pi_t(dr)=E(f(w_1(t))\cdot p(w_1,w_2)/w_2).
$$

In view of the previous lemma it is enough to get the equation for
$\pi_t.$ The next lemma contains the necessary facts from the theory
of extended stochastic integral.

{\bf Lemma 5.2.} {\it Let $H$  be the separable Hilbert space, $\xi$
be a generalized Gaussian random element in $H$ with zero mean and
identity covariation. Suppose, that the random element $x$  in $H$
has two stochastic derivatives and let $I$ and $D$   be the symbols
of the extended stochastic integral and stochastic derivative
correspondingly.  Then for arbitrary $h\in H$  and stochastically
differentiable bounded random variable $\alpha$   the following
formulas hold

1) $\alpha I(x)=I(\alpha x)+(x, D\alpha),$

2) $(DI(x), h)=(x,h)+I((Dx, h)).$ }

{\it Proof.}   The first statement is the well-known relation [12].
Let us check 2).  Use the integration by part formula. Consider the
random variable $\beta$  which is twice stochastically
differentiable. Then, using 1),
$$
E(DI(x),h)\beta =E(DI(x), \beta h)=
$$
$$
= EI(x)I(\beta h) =EI(x)(\beta I(h)-(D\beta, h))=
$$
$$
=E(x, D(\beta I(h)-(D\beta,h)))=
$$
$$
= E[(x, \beta h)+(x, D\beta)I(h)-(x, (D^2\beta, h))]=
$$
$$
= E(x,h)\beta +E(x, I(D\beta h))=
$$
$$
=E\beta((x,h)+I((Dx,h))).
$$
Lemma is proved.

{\bf Remark.}   Note, that the statement of the lemma remains to be
true in the case, when $\alpha$  is not bounded but all terms are
well-defined. Also, due to the formula (4.5) the lemma holds in the
case, when the initial generalized Gaussian random element has not
identity covariation.

Now let us turn  to our filtration problem.

Take the function $f\in C^2_0(\mbR).$  Then for arbitrary $r\in\mbR$
from the Ito formula
$$
f(r+w_1(t))p(w_1,w_2)=
$$
$$
=f(r)p(w_1,w_2)+\int^t_0f'(r+w_1(s))dw_1(s)p(w_1,w_2)+
$$
$$
+\frac{1}{2}\int^t_0f''(r+w_1(s))p(w_1,w_2)ds.
$$
Consider the second summand. It contains the Ito integral which
considers with the extended stochastic integral as it was mentioned
before. So we can apply the  formula 1) from the lemma 5.2:
$$
p(w_1,w_2)\int^t_0f'(r+w_1(s))dw_1(s)=
$$
$$
=\int^t_0f'(r+w_1(s))p(w_1,w_2)dw_1(s)+
$$
$$
+\int^t_0f'(r+w_1(s))(SDp(w_1,w_2))_1(s)ds.
$$
Here the index 1 symbolize the first coordinate of correspondent
element from $H.$  Now note,  that the conditional expectation with
respect to $w_2$  is an operator of the second quantization. So, if
we denote by
$$
\gamma_1(t)=E(w_1(t)/w_2),
$$
then, due to the theorem 1.1
$$
E(f(r+w_1(t))p(w_1,w_2)/w_2)=
$$
$$
=
E(f(r)p(w_1,w_2)/w_2)+\int^t_0E(f'(r+w_1(s))p(w_1,w_2)/w_2)d\gamma(t)+
$$
$$
+\frac{1}{2} \int^t_0E(f''(r+w_1(s))p(w_1,w_2)/w_2)ds+ \eqno(5.2)
$$
$$
+\int^t_0E(f'(r+w_1(s)(SDp(w_1,w_2))_1)(s)/w_2)ds,
$$
where the integral with respect to $\gamma$ is an extended
stochastic integral.
 In order to get the stochastic differentiability of $p$  let us consider
the case when Carleman-Fredholm determinant $\zeta$  is equal to
one. Denote by $P_t$  the orthogonal projector in $L_2([0;1])\oplus
L_2([0;1])$ on the subspace $L_2([0;t])\oplus L_2([0;t]).$

{\bf Lemma 5.3.} {\it Suppose, that the operator $S$  has the
property
$$
\forall t\in[0;1]: \ \ \ \  P_tS=P_tSP_t.
$$
Then $\zeta=1.$ }

{\it Proof.} The value $\zeta$ is the Carleman-Fredholm determinant
of the operator $SDh,$   where
$$
h=(a_1(w_1(\cdot)), a_2(w_1(\cdot))).
$$
Now
$$
Dh(t,s)=
\begin{pmatrix}
a'_1(w_1(t))\1_{[0;t]}(s)&0\\
\\
a'_2(w_1(t))\1_{[0;t]}(s)&0
\end{pmatrix}
\eqno(5.3)
$$

In order to prove that
$$
\det\nolimits_2(Id+SDh)=1
$$
we will use the theorem 3.6.1 from [18].  Due to this theorem it is
enough to check, that the operator $SDh$   is quasi-nilpotent, i.e.,
that
$$
\lim_{n\to\infty}\|(SDh)^n\|^{\frac{1}{n}}=0. \eqno(5.4)
$$
It follows from representation (5.3), that\newline $\forall \vf\in H
\ \forall t\in[0;1]$
$$
\|P_tDh\vf\|^2\leq c\int^t_0\|P_s\vf\|^2ds,
$$
where $c$  depends on the $\sup_\mbR(|a'_1|+|a'_2|).$

Consequently
$$
\|(SDh)^n(\vf)\|^2\leq\|S\|^2\cdot c\cdot
\int^1_0\|P_{t_1}(SDh)^{n-1}(\vf) \|^2dt_1\leq
$$
$$
\leq \|S\|^2\cdot c\cdot\int^1_0 \|P_{t_1}S
P_{t_1}P_{t_1}Dh(SDh)^{n-2}(\vf)\|^2dt_1\leqq
$$
$$
\leq \|S\|^4\cdot c^2\cdot\int^1_0\int^{t_1}_0
\|P_{t_2}(SDh)^{n-2}(\vf)\|^2dt_2dt_1\leq\ldots\leq
$$
$$
\leq \|S\|^{2n}\cdot
c^n\cdot\int^1_0\int^{t_1}_0\ldots\int^{t_{n-1}}_0
\|P_{t_n}\vf\|^2dt_n\ldots dt_1\leq
$$
$$
\leq\|\vf\|^2\frac{\|S\|^{2n}c^n}{n!}.
$$
This means, that (5.4) holds and $\zeta=1.$

Lemma is proved.

Now one can conclude that $p$  has the stochastic derivative and
(5.1)  is correct. The further concretization of (5.2)  can be
possible due to the special form of $p.$   As a  consequence of
(5.2) and (5.4)   we have the following theorem.

{\bf Theorem 5.1.}  {\sl Suppose that the coefficients $a_1, a_2$
and the operator $V$ satisfy the conditions of the lemma 5.3. Then
the random function
$$
U(r,t)=E(f(r+w_1(t))p(w_1,w_2)/w_2)
$$
satisfies relation
$$
dU(r,t)=\frac{1}{2} \frac{\pt^2}{\pt r^2}U(r,t)dt+ \eqno(5.5)
$$
$$
+\frac{\pt}{\pt r} U(r,t)\gamma(dt)+Ef'(r+w_1(t))
(SDp(w_1,w_2))_1(t)dt.
$$
}

In some particular case the last term can be written in a simple
form. For example, when  $a_2=0,$ then (5.5)  transforms into
$$
dU(r,t)=\frac{1}{2} \frac{\pt^2}{\pt r^2} U(r,t)dt+\frac{\pt}{\pt
r}U(r,t)\gamma(dt)+
$$
$$
+a_1(r)\frac{\pt}{\pt r}U(r,t)dt
$$
due to the theorem 3.2.

\vskip 1cm \centerline{REFERENCES}

\begin{enumerate}
\item
It\^o, Kiyosi; McKean, Henry P., Jr. Diffusion processes and their
sample paths. Second printing, corrected. Die Grundlehren der
mathematischen Wissenschaften, Band 125. Springer-Verlag,
Berlin-New York, 1974.
\item
Simon Barry. The $P(\phi)_2$ euclidian (quantum) field theory.
Princeton University Press, 1974.
\item
 Dorogovtsev A.A. Stochastic analysis and random maps in Hilbert
space. VSP Utrecht, The Netherlands, Tokyo, Japan. - 1994. - 110
p.
\item
Dorogovtsev A.A. Anticipating equations and filtration problem.
Theory of stochastic processes, 1997. V.3(19), issue 1 - 2. P. 154
- 163.
\item
Dorogovtsev A.A. Conditional measures for diffusion processes and
anticipating stochastic equations. Theory of stochastic processes,
1998. V.4(20). P. 17 - 24.
\item
 Skorokhod A.V. One generalization of the stochastic integral.
Probability theory and its applications. 1975, V. 20, \#2. P.223 -
237.
\item
 Liptser R. Sh., Shyriaev A.N. Statistics of random processes.
Moskow, Nauka, 1974. 696 p.
\item
Billingsley Patrick. Convergence of probability measures. John
Wiley and Sons, Inc., New York, London, Sydney, Toronto. 1968.
\item
Mandelbrot, B.B., van Ness, J. Fractional Brownian motion,
fractional noises and applications, SIAM Rev., 10:422-437. 1968.
\item
 Tindel, S.; Tudor, C.A.; Viens, F. Stochastic evolution
equations with fractional Brownian motion. Probab. Theory Relat.
Fields 127, No. 2, 186-204.
\item
Dorogovtsev A.A. An action of Gaussian strong random operator on
random elements. Probability theory and its  applications. - 1986.
- V.31,N 4. - P. 811-814.
\item Malliavin, Paul Stochastic analysis. 1997, Text.Mo\-no\-graph,
Grundlehren der Mathematischen Wissenschaften. 313., Berlin:
Springer.
\item
 Nualart, David The Malliavin calculus and related topics. 1995,
Text.Mono\-graph, Probability and Its Applications., New York, NY:
Springer-Verlag.
\item
 Dorogovtsev A.A. Stochastic integration and one class of
Gaussian stochastic processes. Ukr. math. journal, 1998. V.50,\#4.
P.495 - 505.
\item
 Dorogovtsev A.A., Bakun V.V. Random mappings and generalized
additive functionals from Wiener process. Probability theory and
its applications. 2003. V.48, \#1. P.43-61
\item
 Yaozhong Hu. Self-intersection local time of fractional Brownian
motions - via chaos expansion
 J.Math.Kyoto Univ., 41-2 (2001), 233 - 250.
\item
 Buckdahn, R.; Malliavin, P.; Nualart, D. Multidimensional linear
stochastic differential equations in the Skorohod sense. 1997,
Text.Article, Stochastics Stochastics Rep. 62, No.1-2, 117-145
(1997).
\item
A.Suleyman Ustunel, Moshe Zakai.  Transformation of Measure on
Wiener Space. -- Springer, 2000. -- 298 p.
\end{enumerate}

\end{document}